\documentclass[11pt]{article}
    \usepackage[cmex10]{amsmath}\usepackage{amsfonts,amssymb,amsthm}
            \usepackage{xcolor,graphicx,xspace}
    \theoremstyle{plain}       \newtheorem{lemma}{Lemma}
                               \newtheorem{theorem}{Theorem}
                               \newtheorem{corollary}{Corollary}
                               
    \theoremstyle{definition}  \newtheorem{example}{Example}

    \theoremstyle{remark}      \newtheorem{remark}{Remark}
    \DeclareMathAlphabet{\mathsfsl}{OT1}{cmss}{m}{sl}
    \DeclareMathAlphabet{\mathbfsl}{OT1}{cmr}{bx}{it}
    \newfont{\mitq}{cmmi10 at 9pt}    % p265 Kopka1
    \newcommand{\vynech}[1]{}
\newcommand{\field}[1]{\mathbb{#1}}
\newcommand{\R}{\field{R}}

\newcommand{\G}{\mathcal{G}}

\newcommand{\Pd}{{\field{P}}}

\newcommand{\kmax}{{k_{\rm max}}}
\newcommand{\kcr}{{k_{\rm cr}}}

\newcommand{\dom}{{\rm dom}}

\newcommand{\cl}{{\rm cl}}

\newcommand{\thbar}{{\overline{\theta}}}

\newcommand{\intdom}{{\mbox{int\;dom\;}}}
\newcommand{\tet}{{\tilde{\theta}}}
\newcommand{\tete}{{{\theta_1,\theta_2}}}
\newcommand{\tette}{{{\tet_1,\theta_2}}}

\newcommand{\tebteb}{{{\thbar_1,\thbar_2}}}
\newcommand{\ptt}{{p_{\theta_1,\theta_2}}}

\newcommand{\thmin}{{\theta_{\rm min}}}
\newcommand{\qt}{{q_{\theta_2}}}
\newcommand{\qte}{{\hat{q}_k}}
\newcommand{\tetil}{{\tet}_{\rm min}}
\newcommand{\ktil}{{\tilde{k}_{\rm cr}}}

\begin{document}
\title{Almost Worst Case Distributions in Multiple Priors Models\footnote{ICs acknowledges support by the Hungarian National Foundation for Scientific Research OTKA, grant No. K105840. TB acknowledges support by the Josef Ressel Centre for Applied Scientific Computing.}}
\author{Imre Csisz\'{a}r \and Thomas Breuer}
\date{4 June, 2015}

\maketitle

\noindent%

\begin{abstract}
A worst case distribution is %defined to be the minimiser, among the plausible
                             %risk factor distributions, of the expectation of
                             %some random payoff distribution. 
a minimiser of the expectation of some random payoff within a family
 of plausible risk factor distributions.
The plausibility of a risk factor distribution is quantified by a convex 
integral functional. This includes the special cases of relative entropy, 
Bregman distance, and $f$-divergence. 
An ($\epsilon$-$\gamma$)-almost worst case distribution is a risk factor
distribution which violates the plausibility constraint at most by the amount
$\gamma$ and for which the expected payoff is not better than the worst case
by more than $\epsilon$. From a practical point of view the localisation of
almost worst case distributions may be useful for efficient hedging against
them. We prove that the densities of almost worst case distributions cluster
in the Bregman neighbourhood of a specified function, interpreted as worst
case localiser. In regular cases, it coincides with the worst case density,
but when the latter does not exist, 
% But if the infimum of payoff expectations over the plausible distributions
% is not achieved, 
the worst case localiser %is not a worst case density, and 
is perhaps not even a density. %Two other results are an explicit formula for a
                            %computationally cheap
We also discuss the calculation of the worst case localiser, and 
its dependence on the threshold in the plausibility constraint.
 %the characterisation of the worst case expectation also in the pathological cases where the plausibility threshold is so low that it is not achieved by any member of the generalised exponential family.
\end{abstract}

\newpage
\section{Introduction}
Let the monetary payoff or utility of some action, e.g. of a portfolio choice, be described by a function $X(r)$ of a collection $r$ of random risk factors. Suppose the probability distribution which governs the risk factors is not known exactly but may be assumed to belong to a set $\Gamma$ of distributions 
on the sample space $\Omega$ of scenarios $r$ 
(multiple priors model). Then the worst case expected payoff  
\begin{equation}
\label{1}
\inf_{\Pd\in\Gamma} E_\Pd(X) = \inf_{\Pd\in\Gamma} \int_\Omega X(r) \Pd (dr)
\end{equation} 
may be taken as (the negative of) the model risk caused by the lack of knowledge about $\Pd$.
The same expression emerges also in the theory of preferences. Ambiguity averse decison makers may rank possible actions by the criterion of expected utility in the worst case over $\Gamma$. Risk measures or preference criteria of a more general kind involve penalised expected payoff or utility %in the worst case over all probability distributions on $\Omega$,
\begin{equation}
\label{2}
\inf_{\Pd} (E_\Pd(X) + \alpha(\Pd)),
\end{equation} 
where $\alpha(\Pd)$ is a suitable penalty term.  For details, including axiomatic considerations leading to \eqref{1} or \eqref{2}, we refer for example to F\"{o}llmer and Schied \cite{FollmerSchied2004}, 
Hansen and Sargent \cite{HansenSargent2008}, or Gilboa \cite{Gilboa:2009}. 

Any risk measure satisfying some natural postulates (in which case they are dubbed coherent) can be represented as the negative of \eqref{1} for some convex set of distributions $\Gamma$. Relaxing coherence to ``convexity'' yields \eqref{2}, with some convex penalty term $\alpha(\Pd)$. For our purposes, axiomatic theory serves as motivation only. In that theory the infimum in \eqref{1} typically equals a minimum. In models treated in this paper, a worst case distribution $\Pd\in\Gamma$ attaing the minimum in \eqref{1} need not exist. 

If a ``best guess'' $\Pd_0$ of the unknown risk factor distribution is
available, it is natural to use \eqref{1} with $\Gamma$ consisting of those
distributions $\Pd$ that do not deviate much from $\Pd_0$. In the literature
many  measures of deviation of distributions are available; the majority are
non-symmetric. The most versatile one, in various scientific disciplines, is
$I$-divergence or relative entropy. For an axiomatic approach distinguishing 
$I$-divergence in the context of inference 
see Csisz\'{a}r \cite{Csiszar1991} and references therein. Relaxing some 
axioms, that approach leads as alternatives to other frequently used measures
of deviation of distributions, known as $f$-divergences and Bregman distances,
see Section~\ref{sec-Earlier Results} for definitions. In the context of risk
and preferences several authors, perhaps first Hansen and Sargent
\cite{HansenSargent2001}, have considered \eqref{1} with $\Gamma$ equal to an
$I$-divergence ball around $\Pd_0$, or \eqref{2} with $\alpha(\Pd)$ equal to a
constant times the $I$-divergence of $\Pd$ from $\Pd_0$. The preference
relation based on \eqref{2} with this choice of $\alpha(\Pd)$, called
multiplier preferences in \cite{HansenSargent2008}, has been axiomatically
distinguished
 by Strzalecki \cite{Strzalecki2011}. 
Moreover, according to Ahmadi-Javid \cite{Ahmadi-Javid2011} the coherent risk measure he calls entropic value at risk, obtained by taking an $I$-divergence ball for $\Gamma$ in \eqref{1}, is superior to others from the point of view of computability. General $f$-divergences have been employed in this context by Maccheroni {\em et al.}\cite{MaccheroniMarinacciRustichini2006} and Ben-Tal and Teboulle \cite{BenTalTeboulle2007}, see also references in  \cite{BenTalTeboulle2007} to prior work of its authors. Bregman distance could be used similarly but to this we do not have references.

We consider problem \eqref{1} with $\Gamma$ of the following form, including as special cases $I$-divergence balls, $f$-divergence balls and Bregman balls:
\begin{equation}
\label{Gamma1}
\Gamma:= \left\{\Pd: d\Pd=pd\mu,\;  H( p )\leq k\right\}, 
\end{equation}
where $\mu$ is a given measure on $\Omega$ and $H$ is a convex integral functional 
as specified in Section~\ref{sec-Geneneral framework}. A corresponding choice of $\alpha(\Pd)$ in \eqref{2} is 
$\alpha(\Pd) = \lambda H(p)$, $\lambda > 0$. 

Our main focus in this paper is the {\em location} of the infimum, rather than the value of the worst case expected
payoff \eqref{1} or the related infimum \eqref{2}. In cases the infimum is not achieved, there is no worst case distribution, then it is not obvious what the location of infimum should mean. We introduce the concept and prove the existence of a ``localiser of almost worst case distributions'', which in the following sense characterises the location of the infimum, whether or not the minimum is achieved: almost worst case distributions achieving values ever closer to the infimum are in ever smaller Bregman balls around the localiser. Part of the results were presented in the
symposium contribution \cite{BCsISIT}. 

The problem %\eqref{1} 
of minimising $E_\Pd(X)$ subject to $H(p)\leq k$ is related to the problem of
minimising convex integral functionals subject to moment constraints.  This
problem, an extension of 
the celebrated ``information geometric'' problem of $I$-divergence minimisation has been extensively studied in the literature. We rely upon those results in the form presented by Csisz\'{a}r and  Mat\'{u}\v{s} \cite{CsiszarMatus2012art} and we use the basic framework of Breuer and Csisz\'{a}r \cite{BreuerCsiszar2013-MaFi} presented in Section \ref{sec-Earlier Results}.
%Section~\ref{sec-Earlier Results} reviews the basic framework and the results of Breuer and Csisz\'{a}r 
% and  Mat\'{u}\v{s} 
%\cite{BreuerCsiszar2013-MaFi} 
% needed here, are reviewed in Section~\ref{sec-Earlier Results}. 

The new results are presented in Section 3. Theorem 1 in 
Subsection~\ref{sec-Ahmadi} extends a result of Ahmadi-Javid \cite[Theorem
5.1]{Ahmadi-Javid2011} on computing the infimum in \eqref{1} for $\Gamma$
of form \eqref{Gamma1} to our 
%setting, which may be useful for the numerical evaluation of \eqref{1}. 
framework\footnote{This framework does include some assumptions, 
adopted for other purposes. which were absent 
in \cite{Ahmadi-Javid2011}.} that admits also 
non-autonomous integrands and unbounded payoff
function $X$. Our main result, Theorem 2 in 
Subsection~\ref{sec-Almost worst case
  distributions}, addresses
 % contains our main results. The concept of almost worst
                  % case distributions (or densities) is introduced, which
                  % almost achieve the minimum in \eqref{1}. Almost worst case
                  % densities are shown to cluster in the neighbourhood of a
                  % specified function called worst case localiser.  It equals
                  % the density of the worst case distribution attaining the
                  % minimum in \eqref{1}, if that minimum is
                  % attained. Otherwise, the worst case localiser is not a
                  % maximising density and perhaps not a density at
                  % all. 
the worst case and almost worst case distributions (densities) that attain or
almost attain the minimum in \eqref{1}. The almost
worst case densities are shown to cluster, in Bregman distance, around a
specified function called worst case localiser. A similar result is obtained
also for problem~\eqref{2}. The worst case localiser equals the worst case
density if the minimum is attained, while otherwise it is
perhaps not a density at all.
Finally, Subsection~\ref{sec-largek} addresses the effect of the threshold 
$k$ in~\eqref{Gamma1}. Theorems~\ref{WCLWCD} and~\ref{fballs} show that in many situations, including the case of $f$-divergence balls, 
either a worst case distribution exists for all $k>0$ or else it does/does not exist for $k$ less/larger than a critical value 
$\kcr > 0$. It remains open whether a similar result also holds in general---apart from the possibility demonstrated by an example with Bregman balls that no worst case distribution exists for any $k>0$.
%Problem \eqref{1} in case $k$ is larger than a critical value $\kcr$ so that the solution is not a member of the generalised exponential family.

\section{Preliminaries} \label{sec-Earlier Results}
\subsection{General framework}
\label{sec-Geneneral framework}
Let $\Omega$ be any set equipped with a (finite or $\sigma$-finite) measure
$\mu$ on a $\sigma$-algebra not mentioned in the sequel. 
Probability measures $\Pd \ll  \mu$ will be represented by their densities 
$p=d\Pd / d\mu$. The notation $p$ will be used also for nonnegative
(measurable) functions on $\Omega$ which are not densities, i.e., do not have
integral $1$. Equality of functions on $\Omega$ will be meant in the
$\mu$-almost everywhere ($\mu$-a.e.) sense.

Let $H$ be a convex integral functional defined %for (measurable) positive
                                %functions $p$ on $\Omega$ by 
on the vector space of measurable functions\footnote{This functional will be
considered only for nonnegative functions $p$, with no loss of generality
since $p\ge 0\,(\mu$-a.e.)
is a necessary condition for $H( p )<+\infty$, see~\eqref{be0}.}
on $\Omega$ by 
\begin{equation}\label{H1}
H( p ) = H_{\beta,\mu}( p ) := \int_\Omega \beta(r, p( r )) \mu(dr).
\end{equation}
%$p\geq 0$ $\mu$-a.e. is necessary for $H ( p ) < +\infty$, see \eqref{be0}. 
Here $\beta(r,s)$ is a function  of $r\in\Omega,\quad s\in\R$, measurable in 
$r$ for each $s\in\R$, strictly convex and differentiable\footnote{Strict 
convexity appears essential for our main results. Differentiability is assumed
for convenience, it could be dispensed with as in \cite{CsiszarMatus2012art}.}
 in $s$ on $(0,+\infty)$ for each $r\in\Omega$, and satisfying
\begin{equation}
\label{be0}
\beta(r, 0)=\lim_{s\downarrow 0}\beta(r, s), \:\: \:\: \beta(r, s):= +\infty \:\mbox{if $s<0$}.
\end{equation}
Then $\beta$ is a convex normal integrand in the sense of Rockafellar and Wets \cite{RockafellarWets1997}, which ensures the measurability of $\beta(r, p( r ))$ in \eqref{H1} and of similar functions later on. 

Let $X$ be any measurable function interpreted as payoff function, and $\Pd_0$ a default distribution on
$\Omega$ with $\Pd_0 \ll \mu$, $d\Pd_0 / d\mu = p_0$, such that the
expectation 
$$E_{\Pd_0}(X)=\int_\Omega X( r) p_0( r ) \mu(d r ) =: b_0$$ 
exists. Let $m$ and $M$ denote the $\mu$-ess inf and $\mu$-ess sup of $X$, and adopt as standing assumptions 
\begin{eqnarray}
-\infty \leq m & < & b_0 < M \leq +\infty \label{mM}\\
H( p ) & \geq & H(p_0) = 0 \qquad \mbox{whenever} \int p d\mu = 1. \label{H0} 
\end{eqnarray}
%Problem \eqref{1} will be considered for the distribution sets $\Gamma$
%specified in terms of $H$, as in \eqref{Gamma1}.
Due to strict convexity of $\beta$, the inequality in~\eqref{H0} is strict
if $p\neq p_0$.
 %(in the usual $\mu$-a.e. sense).
%Here and later on, equality of functions on $\Omega$ is meant
%$\mu$-almost everywhere.
\begin{example} \label{Example1}
Take $\mu=\Pd_0$, thus $p_0\equiv 1$, and let $\beta(r,s) = f(s)$ be an autonomous convex integrand, with % $f(s)\geq 
$f(1) = 0$ to ensure \eqref{H0}. Then $H( p )$ in \eqref{H1} for 
$d\Pd=pd\mu$ is the {\em $f$-divergence} $D_f(\Pd \, || \,  \Pd_0)$, 
introduced in Csisz\'{a}r \cite{Csiszar1963}. If $f$ is cofinite, i.e. if
$\lim_{s\rightarrow +\infty} f(s)/s = +\infty$, then $\Pd\ll \Pd_0$ is a
necessary condition for  $D_f(\Pd \, || \,  \Pd_0)< +\infty$, hence in that
case $\Gamma$ in 
\eqref{Gamma1} equals the $f$-divergence ball $\{\Pd: D_f(\Pd \, || \,  \Pd_0) \leq k \}$. 
% For $f(s)=\log s+1-s$, $D_f$ equals the Burg entropy.
If $f$ is not cofinite, $f$-divergence may be finite also in 
absence of absolute continuity. Still, with some abuse of terminology, the set
in \eqref{Gamma1} will be called 
$f$-divergence ball also in that case.
%, the $f$-divergence ball may also contain distributions 
%not absolutely continuous to $\Pd_0=\mu$. Then $\Gamma$ consists of the 
%distributions $\Pd  \ll  \Pd_0$ in the $f$-divergence ball.
\end{example}
 
\begin{example} \label{Example2}
Let $f$ be any strictly convex and differentiable function on $(0, +\infty)$, 
and for $s\geq 0$
let $\beta(r, s) = \Delta_f(s, p_0( r ))$ where%\footnote{In~\eqref{delta},
 % $f(0)$ and $f'(0)$ are defined as limits; if
%$f(0)=+\infty$, set $\Delta_f(s,0):=s\times (+\infty).$} 
\begin{equation}
\label{delta}
\Delta_f(s,t):= f(s) - f(t) -f'(t)(s-t).
\end{equation}
Here $f(0)$ and $f'(0)$ are defined as limits; if
$f(0)=+\infty$, we set $\Delta_f(s,0):=0$ for $s=0$ and $\Delta_f(s,0):=\infty$ otherwise.

In this example $\Pd_0 \ll  \mu$ is arbitrary, except that in case $f'(0)=-\infty$ we assume that $p_0 > 0$ $\mu$-a.e.. Then $H( p )$ equals the Bregman distance \cite{Bregman1967}
\begin{equation}
\label{Bf}
B_{f,\mu}(p,p_0):= \int_\Omega \Delta_f(p( r ), p_0( r )) \mu(dr),
\end{equation}
and $\Gamma$ is a Bregman ball of radius $k$ around $\Pd_0$.
Note that here the assumption $f(1)=0$ is not needed to guarantee
\eqref{H0}, but may be adopted anyhow for the function $\Delta_f(s,t)$
is not affected by adding a constant to $f$. 
\end{example}
In the special case $f(s)=s\log s$ %- s + 1
 both examples give the $I$-divergence ball $\Gamma = \{\Pd: D(\Pd \, || \,
 \Pd_0) \leq k \}$ where 
$$ D(\Pd \, || \,\Pd_0):=\int p\log\frac{p}{p_0}d\mu.$$
As another special case, the choice $f(s)=s^2,\,s>0$ gives 
$\Delta_f(s,t)=(s-t)^2$ and $B_{f,\mu}(p,p_0):= \int (p-p_0)^2 d\mu,$ which is the squared $L^2$-distance between $p$ and $p_0$.

\begin{figure}[htbp]
\begin{center}
   \centering
% graphics[width=10.9cm,height=7.5cm]{figures/assetprice-nobond-final-trajectory.pdf}
\includegraphics[height=6cm]{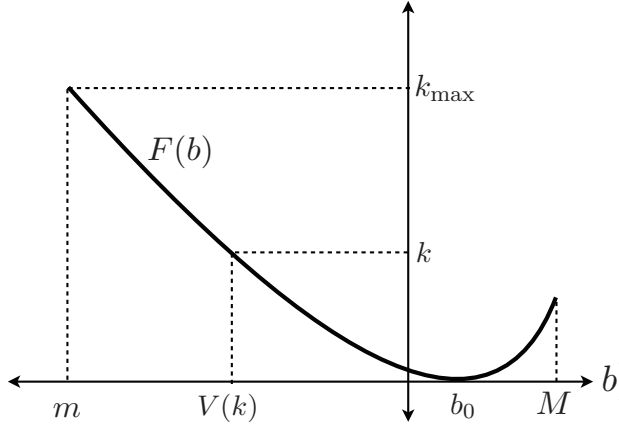}
\end{center}
\caption{\label{fig-F}
Lemma~\ref{Lemma1} relates problem \eqref{V} to the information theoretic problem \eqref{F}: $F(V(k))=k$.
} 
\end{figure}

\bigskip
% \subsection{Basic approach}
For $\Gamma$ of the form \eqref{Gamma1} the infimum in \eqref{1} equals 
\begin{equation}
\label{V}
V( k ) :=  \inf_{p:\int p d\mu =1, H( p ) \leq k} \int X p d\mu,
\end{equation}
and for $\alpha(\Pd):=\lambda H(p)$, $\lambda>0$, the infimum in \eqref{2}
equals
\begin{equation}
\label{W}
W(\lambda):=\inf_{p:\int p d\mu =1}\left[\int Xpd\mu+\lambda H(p)\right].
\end{equation}
The next lemma relates the solution of problem \eqref{V} to that of the following minimisation problem, see Fig.~\ref{fig-F}:
\begin{equation}
\label{F}
F(b):=   \inf_{p:\int p d\mu =1, \int X p d\mu=b} H( p ).
\end{equation}
$F(b)$ is a convex function with minimum $0$ attained at $b=b_0$.
%E_{\Pd_0}(X)=: b_0$. 
% Problems~\eqref{V} and \eqref{F} are closely related. What is the objective
% function in one problem is the constraint in the other, and vice versa (see
% Fig.~\ref{fig-duality}).  The solutions of the two problems are related in
% the following way: 
A standing assumption will be, in addition to \eqref{mM},\eqref{H0}, that
\begin{equation}
\label{H2}
\kmax:= \lim_{b\downarrow m} F(b)>0.
\end{equation}
This is a necessary condition for the functional $H$ to yield a nontrivial
measure of risk for the payoff function $X$, since $\kmax=0$ would imply
$V(k)=m$ for each $k>0$. Note that if $m=-\infty$ then $\kmax=+\infty$
(subject to \eqref{H2}), while if $m$ is finite then $\kmax\le F(m)$ where
the strict inequality is possible.

\begin{lemma}{\cite[Proposition 3.1]{BreuerCsiszar2013-MaFi}} 
\label{Lemma1}
%If $0 < k < \kmax:= \lim_{b\downarrow m} F(b)$, 
To each $k\in(0,\,\kmax)$ there exists a unique $b \in (m,b_0)$ %satisfying 
%\begin{equation}\label{Fb}
%F(b)=k, % \qquad m < b < b_0
%\end{equation}
with $F(b)=k$,  and then $V( k ) = b.$ %\qquad 0 < k < \kmax$$
%\end{lemma} %For a proof see Breuer and Csisz\'{a}r . 
The minimum in \eqref{V} is attained if and only if that in \eqref{F} is attained (for the above $b$), % in \eqref{Fb}
and then the same $p$ attains both minima.
\end{lemma}

\begin{remark}\label{Remark1}
The %condition $0 < k < \kmax$ 
assumption on $k$ is not restrictive, for if $k=0$ or $ k\geq \kmax >0$ then 
$V(k)$ trivially equals $b_0$ or $m$.
%trivially $V(k)=b_0$ resp. $V(k)=m$. If $\kmax=0$ then the functional $H$ is unsuitable to define a non-trivial measure of risk for the given $X$, as then $V(k)=m$ for each $k>0$. 
\end{remark}
\begin{remark}\label{Remark2}
By \cite[Theorem 2]{BreuerCsiszar2013-MaFi}, the standing assumption \eqref{H2}
is equivalent to \eqref{nont} below (which automatically holds if
$m>-\infty$),
and that condition implies $F(b)>0$
for each $b<b_0$. In particular, the continuous convex function 
$F(b),\;b\in(m,b_0]$ is strictly decreasing, and $V(k), \;k\in [0,\kmax)$ 
equals its inverse function. 
%The latter implies, in turn, that the function $V(k)$ is
%right continuous at $k=0$, and the hypothesis of Lemma~\ref{Lemma1} could be 
%weakened
%to $k\in[0,\,\kmax)$ if the interval $(m,b_0)$ in the assertion were replaced
%by $(m,b_0].$ 
\end{remark}
 
\subsection{Basic concepts and facts}
Lemma~\ref{Lemma1} admits to treat Problem~\eqref{V} using known results about minimising convex integral functionals under moment constraints, specifically with moment mapping defined by $\phi( r ):= (1, X( r ))$. We will rely upon results in Csisz\'{a}r and  Mat\'{u}\v{s} \cite{CsiszarMatus2012art}\footnote{Many of these results have been known earlier, though typically under less general conditions.}, specified for this moment mapping. Then the value function in \cite{CsiszarMatus2012art} becomes
\begin{equation}
\label{J}
J(a, b) :=  \inf_{p:\int p d\mu =a, \int X p d\mu =b} H( p ), 
\end{equation}
thus $F(b)=J(1,b)$. 

%\subsection{Basic concepts and facts}
%Some basic facts needed from Breuer and Csisz\'{a}r
%\cite{BreuerCsiszar2013-MaFi} and Csiszar and Mat\'{u}\v{s}
%\cite{CsiszarMatus2012art} are recalled, and some consequences are added. 
The function $J$ in \eqref{J} is convex, and its effective domain $\dom \, J:=\{(a,b): J(a,b)< +\infty\}$ has interior
\begin{equation}
\label{eq-intdomJ}
\intdom J = \{(a,b): a>0, \, am < b < aM\},
\end{equation}
by \cite[Lemma 6.6]{CsiszarMatus2012art}. The function $J$ is proper (not identically $+\infty$ and never equal to $-\infty$) because it equals zero at $(1,b_0)\in \intdom J $, see \eqref{mM}, \eqref{H0}.
Hence its convex conjugate $J^*(\theta_1, \theta_2):= \sup_{a,b}[\theta_1 a+
\theta_2 b -J(a,b)]$ is a closed (i.e., lower semicontinuous) proper
convex function. A crucial fact is the instance of 
\cite[Theorem 1.1]{CsiszarMatus2012art} that
\begin{equation}
\label{K}
J^*(\tete) = K(\theta_1, \theta_2) :=  \int \beta^*(r, \theta_1 + \theta_2 X( r ))\mu(dr),  
\end{equation}
where $\beta^*$ is the convex conjugate of $\beta$,
\begin{equation}
\label{bstar}
\beta^*(r,\tau):= \sup_{s\in \R}\left(s\tau - \beta(r,s)\right).
\end{equation} 
The conjugate and derivatives of $\beta$ are by the second variable. 

Below, derivatives at $0$ and $+\infty$ are interpreted as limits of derivatives at $s\downarrow 0$ and $s \uparrow +\infty$.
For fixed $r\in\Omega$ the function $\beta^*$ equals $-\beta(r,0)$ for $\tau\leq \beta'(r,0)$, it is strictly convex in the interval 
$(\beta'(r,0), \beta'(r,+\infty))$, and equals $+\infty$ if $\beta'(r,
+\infty)$ is finite and $\tau>\beta'(r,+\infty)$. 
This function is differentiable in the interval $(-\infty, \beta'(r,+\infty))$.
Its dervative $(\beta^*)'(r,\tau)$ is positive
 and strictly increasing in $(\beta'(r,0),\beta'(r,+\infty))$, and
approaches $0$ or $+\infty$ as $\tau\downarrow 
\beta'(r,0)$ or $\tau\uparrow \beta'(r,+\infty)$.

Since $J^*=K$ implies $J^{**}=K^*$, and $J^{**}$ (equal to the closure of $J$) 
may differ from $J$ only on the boundary of $\dom\, J$,
\begin{equation}
\label{FK}
F(b)=J(1,b) = K^*(1,b) = \sup_{\tete} [\theta_1+\theta_2 b - K(\tete)],
\end{equation}
except possibly for $b$ equal to $m$ or $M$, see \eqref{eq-intdomJ}.
%Moreover, if $m < b < M$ then in \eqref{FK} the maximum is attained:  
%$(\tete)$ attains it if and only if the subdifferential\footnote{The 
%subdifferential $\partial f(x)$ of a proper convex function 
%$f$ at a point $x$ is the set of subgradients of $f$ at $x$, i.e. the set of 
%vectors $x^*$ such that $f(z) \ge f(x) + \langle x^*,z-x\rangle$ for all $z$. 
%If $f$ is differentiable at $x$, its unique subgradient at $x$ is the gradient 
%$\nabla f(x)$. The differentiability of the function $J$ is not addressed here.%}
%%proper and $x\in \intdom f$, then $\partial f(x)\not= \emptyset$, see 
%%\cite[Theorem 23.4]{Rockafellar1970}.} 
%$\partial J (1,b)$ contains $(\tete)$, or equivalently, $\partial K(\tete)$
%contains $(1,b)$, see \cite[Theorem 23.5]{Rockafellar1970}; the
%former subdifferential is nonempty by \cite[Theorem 23.4]{Rockafellar1970}.
%Note that \eqref{FK} 
This can be rewritten as
\begin{equation}
\label{FL}
F(b) = \sup_{\theta_2}\;[\theta_2 b - G(\theta_2)] = G^*(b)
\end{equation}
where
\begin{equation}
\label{L}
G(\theta_2) := \inf_{\theta_1}\;[K(\tete)-\theta_1].
\end{equation}
The function $G$ will play a similar role as the logarithmic moment generating
function does when $\Gamma$ in~\eqref{Gamma1} is an 
$I$-divergence ball, see Example~\ref{Kullback}. A consequence of~\eqref{FL} applied to $b=b_0$: is the simple bound
\begin{equation}\label{Gbound}
  G(\theta_2)\ge \theta_2 b_0.
\end{equation}

The following family of non-negative functions on $\Omega$ will play a key role like exponential 
families do for $I$-divergence minimisation:
\begin{equation}\label{EF}
p_{\theta_1, \theta_2}( r ):=  (\beta^*)'(r, \theta_1 + \theta_2 X( r )), \:\:\: (\theta_1, \theta_2)\in\Theta
\end{equation}
where%\footnote{The definition~\eqref{Te} makes sure that the derivative in \eqref{pte} exists for $\mu$-a.e. $\omega\in\Omega$ if $(\theta_1, \theta_2)\in\Theta$. For all other $\omega\in\Omega$, if any, one may set by definition $p_{\theta_1, \theta_2}=0$.} 
\begin{equation}
\label{Th}
\Theta:= \left\{ (\theta_1, \theta_2)\in\dom \; K: \theta_1 + \theta_2 X(r ) < \beta'(r, +\infty) \quad \mu\mbox{-a.e.}\right\}.
\end{equation}
% with $\dom\;K:= \{(\theta_1, \theta_2): K(\theta_1, \theta_2) < +\infty\}.$

\begin{remark} \label{unique} It may happen that different parameters 
$(\tete)\in\Theta$ give rise to equal functions \eqref{EF}, but only in case of 
functions that equal zero except for $r$ in a set where $X(r)$ is 
constant $\mu$-a.e. This follows because for any fixed $r\in\Omega$,  
the fact that $(\beta^*)'(r,\tau)$ is strictly increasing for
$\tau\in (\beta'(r,0),\beta'(r,+\infty))$ implies that $p_\tete (r)$
in \eqref{EF}, if positive, uniquely determines $\theta_1 +\theta_2 X(r)$. 
In particular, for positive valued functions \eqref{EF} the parameters 
$(\tete)\in\Theta$ are always unique, due to the standing assumption \eqref{mM}.   
\end{remark}

As $(\tete)\in\dom\, K$ implies $(\tette)\in\Theta$ for each
$\tet_1<\theta_1$, the sets $\dom\;K$ and $\Theta$ have the same projection
to the $\theta_2$-axis. This projection will be denoted by $\Theta_2$.
It is a (finite or infinite) interval. The standing
assumptions~\eqref{mM},~\eqref{H0} imply that $\Theta_2$
contains the origin, and the default density
%some negative numbers, too. Indeed, $0\in\Theta_2$ follows from the standing
%assumptions~\eqref{mM},~\eqref{H0}, which also imply that the default density 
$p_0$ belongs to the family~\eqref{EF} with $\theta_2=0$, 
see \cite[Remark 4]{BreuerCsiszar2013-MaFi}. The left endpoint of the interval 
$\Theta_2$ will be denoted by $\theta_{\min}$. By 
\cite[Theorem 2]{BreuerCsiszar2013-MaFi}, the standing assumption $\kmax>0$
is equivalent to 
\begin{equation}\label{nont}
  % \Theta_2 \quad\mbox{contains some}\quad \theta_2<0
  \theta_{\min}<0.
\end{equation}
%is equivalent to the standing assumption $\kmax>0$, by 
%\cite[Theorem 2]{BreuerCsiszar2013-MaFi}.
%Clearly, $\Theta_2$ is an interval, and it contains the origin, see
%\cite[Remark 4]{BreuerCsiszar2013-MaFi}, the default density $p_0$ being
%a member of the family~\eqref{EF} with $\theta_2=0$. Moreover, by
%\cite[Theorem 2]{BreuerCsiszar2013-MaFi}, our standing assumption $\kmax>0$
%is equivalent to the condition
%\begin{equation}\label{nont}
%   \Theta_2 \quad\mbox{contains some}\quad \theta_2<0.
%\end{equation}
% [Further, $\ptt(r)$ is a non-decreasing function of the parameter $\theta_1$, for each $\theta_2$ and $r\in\Omega$.]

By \cite[Lemma 3.6]{CsiszarMatus2012art}, the directional derivatives of the 
function $K$ in~\eqref{K} can be expressed,
at any $(\tete)\in\Theta$ and
%and $\tet=(\tet_1,\tet_2)\in \dom\, K$
for any $(\tet_1, \tet_2)\in\dom\, K$, as
%in terms of functions in~\eqref{EF}. 
\begin{eqnarray}\label{dide}
%$$K'(\theta,\tet-\theta) := 
\lim_{t\downarrow 0}\frac{1}{t} \left[ K(\theta_1+t(\tet_1-\theta_1),
\theta_2 + t(\tet_2-\theta_2))-K(\tete)\right] \nonumber \\
=\int \left[\tet_1 - \theta_1 + (\tet_2 - \theta_2) X(r)\right]\ptt(r)\mu(dr),
\end{eqnarray}
where the integral is well-defined and is not equal to $+\infty$.
In particular, $K$ is differentiable in the interior of its effective domain,  
with %\cite[Corollary 3.8]{CsiszarMatus2012art} 
%for $(\tete)\in\Theta$ 
\begin{eqnarray}
%\label{grad}
\frac{\partial}{\partial\theta_1} K(\tete)  & = & \int p_{\theta_1, \theta_2} d\mu, \label{par1}\\
\frac{\partial}{\partial\theta_2} K(\tete) & = &  \int X p_{\theta_1, \theta_2} d\mu. \label{par2}
\end{eqnarray}
%\cite[Corollary 3.8]{CsiszarMatus2012art}.
The same equations hold at $(\tete)\in\Theta$ on the boundary
%At boundary points $(\tete)$ 
of $\dom\, K$ 
%that belong to $\Theta$, the same equations hold 
for those one-sided partial derivatives of $K$ which are 
defined there, thus \eqref{par1} holds for the left partial
derivative at each $(\tete)\in\Theta$.

The following lemma gives relevant information about %the minimization 
%problem~\eqref{L}.
evaluating the function $G$ in~\eqref{L}.  
Its proof is effectively contained in the proof of
\cite[Proposition 2]{BreuerCsiszar2013-MaFi}, but for convenience a full
proof will be given in the Appendix. 
Clearly, $\dom\, G:=\{\theta_2: G(\theta_2) < +\infty\} = \Theta_2.$
%If $(\tebteb)$ attains the maximum in \eqref{FK} (since $(\tebteb)\in \partial J(1,b)$),
%then $\thbar_2$ attains the maximum in \eqref{FL} and $\thbar_1$ attains the minimum in \eqref{L} for $\thbar_2$. The next lemma, proved in the Appendix, gives more information about the function $L$.
\begin{lemma}
\label{lemma-L}
Given any $\theta_2\in \Theta_2$, either (i) some
 $\bar{\theta}_1\in \R$ satisfies $(\bar{\theta}_1,\theta_2)\in\Theta$, 
$\int p_{\bar{\theta}_1,\theta_2} d\mu = 1$, %see \eqref{par1}. Then this
                                %$\theta_1$ attains the minimum in \eqref{L}
                                %and $L(\theta_2) = K(\tete)-\theta_1$. 
or (ii) $\tet_1:=\sup \{\theta_1: (\tete)\in \dom\, K\}$ is finite,
 $(\tet_1, \theta_2)\in \Theta$, $\int p_{\tet_1, \theta_2} d\mu < 1$. In 
either case, in \eqref{L} the minimum is attained, and the unique
minimiser is $\bar{\theta}_1$ respectively $\tet_1$.
%attains the minimum in \eqref{L} and $L(\theta_2) = K(\tet_1, \theta_2)-\tet_1$.
\end{lemma}
%%%%%%%%%%%%%
\subsection{Generalised Pythagorean identity}

Given the convex integrand $\beta$, define
$\Delta_{\beta(r, \cdot)}(s,t)$ as %$\Delta_f (s,t)$ in
in \eqref{delta}, 
with the convex function $\beta(r, \cdot): s\mapsto \beta(r, s)$
playing the role of $f$. The mapping $(r, s,t)\mapsto \Delta_{\beta(r,
  \cdot)}(s,t)$ is a normal integrand \cite[Lemma 2.10]{CsiszarMatus2012art},
hence if $p$ and $q$ are non-negative measurable functions on $\Omega$ 
then so is also $\Delta_{\beta(r, \cdot)}(p( r ), q( r ))$, denoted briefly
by $\Delta_\beta(p,q)$. Extending the concept of Bregman distance \eqref{Bf},
% the Bregman distance $B_\beta (p,p_0)$ of measurable, non-negative functions
% $p$, $p_0$ on $\Omega$ 
define 
\begin{equation}
\label{Bbe}
B(p,q)=   B_{\beta,\mu} (p,q) :=  \int \Delta_\beta(p,q) d\mu.
\end{equation}
Like its special case in \eqref{Bf}, it is non-negative and equals $0$ only if
$p=q$. If $\beta=\Delta_f$, as in Example~\ref{Example2}, then $B_{\beta,\mu}$ is equal to the $B_{f,\mu}$ of \eqref{Bf}.

The following lemma, crucial for this paper, is an instance of 
\cite[Lemma 4.15]{CsiszarMatus2012art}, combined with 
\cite[Remark 4.13]{CsiszarMatus2012art}
\begin{lemma}\label{Pyti}
For each density $p$ with $\int Xp d\mu $ finite, and each $(\tete)\in\Theta$,
\begin{eqnarray}\label{Py}
H(p)=\theta_1+\theta_2\int Xp d\mu -K(\tete)+B(p,p_\tete)\nonumber\\
   +\int |\beta'(r,0)-\theta_1-\theta_2 X(r)|_+p(r) \mu(dr).
\end{eqnarray}
\end{lemma}

If $p_\tete$ is a density, the special case $p=p_\tete$ of~\eqref{Py} (or
direct calculation) gives that
\begin{equation}\label{Pyt1}
H(p_\tete)=\theta_1+\theta_2 \int X p_\tete d\mu -K(\tete).
\end{equation}
Then \eqref{Py} and \eqref{Pyt1} imply
\begin{equation}\label{Pyt2}
H(p)=H(p_\tete)+B(p,p_\tete)+\int
|\beta'(r,0)-\theta_1-\theta_2 X(r)|_+p(r) \mu(dr)
\end{equation}
for each density $p$ satisfying
\begin{equation}\label{cond}
\int Xpd\mu=\int Xp_\tete d\mu.
\end{equation}
%As a consequence, if $(\tete)\in\Theta$ is such thar
%\begin{equation}\label{opcond1}
%\int p_\tete d\mu=1,\quad \int Xp_\tete d\mu=b
%\end{equation}
%then $p=p_\tete$ attains the minimum in the definition~\eqref{F} of $F(b).

Identities like \eqref{Pyt2} frequently occur in the literature, primarily
in cases when the last term vanishes (it trivially does if 
$\beta'(r,0)=-\infty$). They are referred to as Pythagorean 
identities,\footnote{If $\beta(r,s)=f(s)=s^2-1\; (s>0)$ then
  $p_\tete=\frac{1}{2}|\theta_1+\theta_2 X(r)|_+$ and \eqref{Pyt2} reduces to
the classical Pythagorean identity $||p||^2=||p_\tete||^2+||p-p_\tete||^2$
provided that \eqref{cond} holds for $\tete$ with 
$\theta_1+\theta_2 X(r)\ge 0.$} and \eqref{Py} will be called generalised
Pythagorean identity.

The above results admit a short proof of the following key lemma, see
\cite[Theorem 1]{BreuerCsiszar2013-MaFi} for a related result.
\begin{lemma}\label{l-opcond}
(i) Let $(\tete)\in\Theta$, $\int p_\tete d\mu =1$. Then 
$\int X p_\tete d\mu$ is finite if and only if $H(p_\tete)$ is. In that case
%and $H(p_\tete)$ finite, 
the density $p=p_\tete$ uniquely\footnote{Uniqueness is meant for the
function, in the $\mu$-a.e. sense. See Remark~\ref{unique}.}
attains the minimum in the
definition~\eqref{F} of $F(b)$ for $b:=\int Xp_\tete d\mu$, and
\begin{equation}\label{Fb=}
F(b)=H(p_\tete)=\theta_1+\theta_2 b-K(\tete)=\theta_2 b-G(\theta_2).
\end{equation}
Supposing
%\footnote{It is left open whether $p_\tete$ with
%$\theta_2\neq 0$ might equal the default density $p_0$ at all, and more 
%generally, 
%whether different parameter pairs $\tete$ might ever yield the
%same density $p_\tete$, in the $\mu$-a.e. sense. The uniqueness assertion of
%the lemma means uniqueness of the function rather than of the parameters.}
 $H(p_\tete)>0$,  here $b$ is less or larger than $b_0$ according as
$\theta_2$ is negative or positive.

(ii) For $k\in (0,\kmax)$, a density $p$ attains the minimum in the definition
\eqref{V} of $V(k)$ if and only if $p=p_\tete$ for some $(\tete)\in\Theta$
with $\theta_2<0$ and
\begin{equation}\label{H=k}
 H(p_\tete)=k\quad\mbox{or equivalently}\quad \int Xp_\tete d\mu =V(k).
%\theta_1+\theta_2\int Xp_\tete d\mu -K(\tete)=k.
\end{equation}
%In particular, if there exists $(\tete)$ as above then
%$V(k)=\int Xp_\tete d\mu$.
\end{lemma}
\begin{proof} The first assertion holds by~\eqref{Pyt1}, and the second one 
since~\eqref{Pyt2},~\eqref{cond} imply $ H(p)> H(p_\tete)$ for
each density $p$ with
$\int Xp d\mu=b=\int Xp_\tete d\mu$. %satisfies $ H(p)> H(p_\tete)$. This proves
%the first assertion, and
Then~\eqref{Fb=} follows by~\eqref{Pyt1} and the
consequence $K(\tete)-\theta_1=G(\theta_2)$ of Lemma \ref{lemma-L}.
Finally, \eqref{Gbound}
%~\eqref{FL} applied to $b_0$ in the role of $b$ gives
%\begin{equation}\label{0ineq} 
%0=F(b_0)\ge \theta_2 b_0-G(\theta_2).
%\end{equation}
and~\eqref{Fb=} yield $\theta_2 b >\theta_2 b_0$, proving the last
assertion of part (i).

(ii) For sufficiency, it is enough to verify the equivalence~\eqref{H=k},
under the given hypotheses. The function $V: (0,\kmax)\rightarrow (m,b_0)$
is the inverse of $F: (m,b_0)\rightarrow (0,\kmax)$, see Lemma~\ref{Lemma1}
and Remark~\ref{Remark2}. This and the result $F(\int Xp_\tete d\mu)=H(p_\tete)$
of part (i) imply~\eqref{H=k}, because if $0<H(p_\tete)<\kmax$ then 
$m<\int Xp_\tete d\mu<b_0$ (the upper bound follows from $\theta_2<0$, due to
the last assertion of part (i)).
%(ii) Since $k\in (0,\kmax)$ and $\theta_2<0$, 
%the condition $b\in (m,b_0)$ in Lemma~\ref{Lemma1} holds for 
%$b=\int Xp_\tete d\mu$. Then part (i) and Lemma~\ref{Lemma1} imply the
%sufficiency assertion, including the equivalence in~\eqref{H=k} under the 
%conditions in~\eqref{H=k} are equivalent by~\eqref{Pyt1}.
%By part (i) and Lemma~\ref{Lemma1}, if $p_\tete$ is a density
%with $H(p_\tete)=k$ then $V(k)=b:=\int Xp_\tete d\mu$ and
%$p=p_\tete$ attains the minimum in~\eqref{V}, 
%provided that the condition $b\in (m,b_0)$ of Lemma~\ref{Lemma1} is met.    
%Since $F(b)=\theta_2 b-G(\theta_2)$ by~\eqref{Fb=}, and $F(b')\ge\theta_2
%b'-G(\theta_2)$ for each $b'\in (m,M)$ by~\eqref{FL}, the hypothesis 
%$\theta_2<0$ implies $b<b'$ for each $b'$ with $F(b')<F(b)$, in particular,
%$b<b_0$. As $b>m$ follows from $H(p_\tete)=k<\kmax$, the
%sufficiency assertion of (ii) is proved. 
Regarding necessity, a
density $p$ that attains the minimum in \eqref{V} clearly satisfies
the constraint $H(p)\le k$ with the equality. We skip the proof of the 
remaining 
assertion that $p$ has to be of form $p_\tete$ with $\theta_2<0$,
for this will be an immediate consequence of Theorem~\ref{main}.
\end{proof}

\section{New results}
\subsection{Calculating V(k)}\label{sec-Ahmadi}

A procedure to calculate $V(k)$ in~\eqref{V} is to first determine the function
$K$ in~\eqref{K}, then the function $F$ via~\eqref{FK} (this may be done in
two steps, first determining the function $G$ in~\eqref{L}), and finally
$V(k)$ as the solution $b\in(m,b_0)$ of the equation $F(b)=k$, see Lemma 1.
%The calculation of $W(\lambda)$ in~\eqref{W} is somewhat less costly, it
%requires the calculation of $G(\theta_2)$ only for a single argument 
%$\theta_2$.
%Indeed, for $\lambda>0$,
%\begin{eqnarray}\label{W=}
%W(\lambda)=\inf_b [b+\lambda F(b)]=-\lambda\sup_b [-\frac{b}{\lambda}-F(b)]
%\nonumber\\
%          =-\lambda F^*(-\frac{1}{\lambda})=-\lambda G(-\frac{1}{\lambda}).
%\end{eqnarray}
%The last equality follows since $F^*=G$, as suggested by~\eqref{FL}. Note
%that~\eqref{FL} would admit $F^*(\theta_2)\neq G(\theta_2)$ when $\theta_2$
%is an endpoint of $\dom\;G =\Theta_2$. A full proof of $F^*=G$ will be 
%given in the Appendix. 
In regular cases, $b=V(k)$ is characterised by equations involving partial
derivatives of the function $K$, see \cite[Corollary 1]{BreuerCsiszar2013-MaFi},
which may facilitate its computation. The following Theorem, combined with 
Lemma~\ref{lemma-L}, may help to reduce computational complexity even in
``irregular''cases. Previously, Ahmadi-Javid \cite[Theorem 5.1]{Ahmadi-Javid2011}
proved an identity equivalent to~\eqref{mm} for autonomous integrands and 
bounded payoff functions.
%\footnote{Our framework includes assumptions 
%about the integrand not made there, which we need for 
%other purposes and appear dispensible for the proof of~\eqref{mm}.}
 
A lemma is sent forward that will be proved in the Appendix.
\begin{lemma}\label{Fstar} $F^*=G$.
\end{lemma}
Note that while 
$F^*=G^{**}=\cl{G}$ immediately follows from~\eqref{FL}, it appears nontrivial 
that the function $G$ is closed. 

\begin{theorem}\label{AJ}
For $k\in(0,\kmax)$
\begin{equation}\label{mm}
V(k)=\max_{\theta_2<0}\max _{\theta_1\in\R}\frac{k+K(\tete)-\theta_1}{\theta_2}
  = \max_{\theta_2<0}  \frac{k+G(\theta_2)}{\theta_2}.
\end{equation}
A maximiser for the second maximum in~\eqref{mm}
is equivalently a maximiser of $\theta_2 b-G(\theta_2)$ where $b=V(k)$.
A pair $(\tete)$ attains the first maximum in~\eqref{mm} if and only
if it attains the maximum in~\eqref{FK}, for
$b=V(k)$. Such $(\tete)$ belongs to $\Theta$ and 
satisfies $\int p_\tete d\mu \le 1$.
%a pair $(\tete)$ 
%$(\tete)\in\partial J(1, V(k))$, or equivalently
%$(1, V(k))\in\partial K(\tete)$. Moreover, this maximiser $(\tete)$
%$\int p_\tete d\mu \le 1$. A maximiser for the second maximum in~\eqref{mm}
%is equivalently a maximiser of $\theta_2 b-G(\theta_2)$ where $b=V(k)$.
\end{theorem}
%The intuitive meaning of the function $p_\tete$ will be
%clarified in Subsection ?
\begin{proof} The conditions $b=V(k)$,  
$k\in(0,\kmax)$ are equivalent to $F(b)=k$,
$b\in(m,b_0)$, see Lemma~\ref{Lemma1} and Remark~\ref{Remark2}. 
The condition that $\theta_2$ is a
maximiser of $\theta_2 b-G(\theta_2)$ means, by~\eqref{FL}, that
\begin{equation}\label{=k}
\theta_2 b-G(\theta_2)=G^*(b)=F(b)=k
\end{equation} 
or equivalently, see Lemma~\ref{Fstar}, that $\theta_2 b-F(b)=F^*(\theta_2)$.
This proves that $\theta_2$ is a
maximiser of $\theta_2 b-G(\theta_2)$ if and only if\footnote{Here $F_{-}'$ and
$F_{+}'$ denote one-sided derivatives; the differentiability of the
function $F$ is not addressed.}
%\begin{equation}\label{derF}
$F'_{-}(b)\le\theta_2\le F'_{+}(b)$.
%\end{equation}
In particular, a (perhaps non-unique) maximiser $\theta_2<0$ does exist.

By~\eqref{=k}, the maximum of $\theta_2 b-G(\theta_2)$ equals $k$, hence
$\theta_2 b-G(\theta_2)\le k$ for each $\theta_2\in\R$. This proves the assertions 
that the maximum of
\begin{equation}\label{slope}
\frac{k+G(\theta_2)}{\theta_2}\quad\quad (\theta_2<0)
\end{equation}
is equal to $b=V(k)$, and a maximiser of~\eqref{slope} is equivalently a 
maximiser of $\theta_2 b-G(\theta_2)$. The remaining assertions 
of Theorem~\ref{AJ}
immediately follow from this and Lemma~\ref{lemma-L}.
\end{proof}

%By Lemma~\ref{Fstar}, $G$ is a closed convex function with
%$G(0)=0$. Let
%\begin{equation}\label{graph}
%\G :=\{(\theta_2, G(\theta_2)): \theta_2\in\Theta_2,\:\theta_2<0\}
%\end{equation}
%denote the graph of $G$ restricted to negative arguments. 
%As
%\begin{equation}\label{slope}
%\frac{k+G(\theta_2}{\theta_2}\quad (\theta_2<0)
%\end{equation}
%equals the slope of of the straigh line through $(0,-k)$ and
%$(\theta_2, G(\theta_2))\in\G$, its maximum equals the slope of
%the supporting line to $\G$ through $(0,-k)$. Denote this maximum
%by $b$, then
%\begin{equation}\label{supp}
%k=G^*(b)=F(b),
%\end{equation} 
%the first equality by the definition of convex conjugate and the second
%one by~\eqref{FL}. 

%Since $G(\theta_2)\ge \theta_2 b_0$, see~\eqref{0ineq}, the slope~\eqref{slope}
%is always less than $b_0$, in particular, $b<b_0$. Hence
%Lemma~\ref{Lemma1} and~\eqref{supp} imply that $b=V(k)$. This proves
%the second equality in~\eqref{mm}, and the first
%one follows since $G(\theta_2)=\max_{\theta_1\in\R}[K(\tete)-\theta_1]$ by  
%Lemma~\ref{lemma-L}. 

%Of the remaining assertions of the Theorem, it suffices to prove the 
%last one, the others follow by Lemma~\ref{lemma-L}. Now, the condition
%on $\theta_2$ that it attains the maximum $V(k)$ of~\eqref{slope} can be
%equivalently written, denoting $V(k)=b$ thus $k=F(b)$, as
%$$ b\theta_2-G(\theta_2)=k=F(b).$$
%Since $F(b)=G^*(b)$, see~\eqref{FL}, this means exactly that $\theta_2$
%%is a maximiser of $b\theta_2-G(\theta_2$. 
%\end{proof}

The calculation of $W(\lambda)$ in~\eqref{W} is somewhat less costly
than that of $V(k)$. It
requires the calculation of $G(\theta_2)$ only for a single value of
$\theta_2$, since for $\lambda>0$ we have (using Lemma~\ref{Fstar} in the
final step)
\begin{eqnarray}\label{W=}
W(\lambda)=\inf_b [b+\lambda F(b)]=-\lambda\sup_b [-\frac{b}{\lambda}-F(b)]
\nonumber\\
          =-\lambda F^*(-\frac{1}{\lambda})=-\lambda G(-\frac{1}{\lambda}).
\end{eqnarray}

\begin{figure}[htbp]
\begin{center}
   \centering
% graphics[width=10.9cm,height=7.5cm]{figures/assetprice-nobond-final-trajectory.pdf}
\includegraphics[height=6cm]{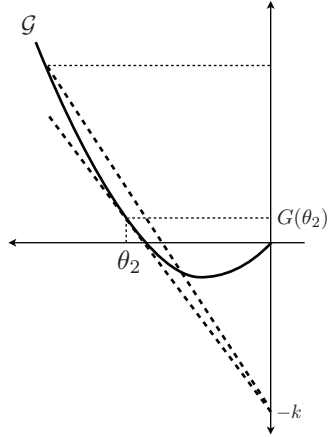}
\end{center}
\caption{\label{fig-G}
The supporting line has maximum slope $b=(k+G(\theta_2))/\theta_2$ among all lines passing through $(0,-k)$ and some point of $\G$. This slope is equal to the solution $V(k)$ of problem~\eqref{V}.
} 
\end{figure}

\begin{remark}\label{geom}
The following geometric interpretation of the proof of Theorem~\ref{AJ}
deserves emphasis, see Fig.~\ref{fig-G}. Denote by
\begin{equation}\label{graph}
\G :=\{(\theta_2, G(\theta_2)): \theta_2\in\Theta_2,\:\theta_2\le 0\}
\end{equation}
the graph of $G$ restricted to nonpositive arguments. Recall that, by
Lemma~\ref{Fstar}, $G$ is a closed convex function with $G(0)=0$. Then
\eqref{slope} is the slope of the straight line through $(0,-k)$ and
$(\theta_2, G(\theta_2))\in\G$, which is maximised by
the supporting line to $\G$ through $(0,-k)$. The proof of Theorem~\ref{AJ}
shows that this supporting line (exists and) has slope $b=V(k)$. The
maximum $b=V(k)$ of~\eqref{slope} is attained if and only if
$(\theta_2, G(\theta_2))$ is on this supporting line. 

%Note for later reference that, by the same proof, also for $k=0$
%the supremum of~\eqref{slope} equals the slope of a 
%supporting line to $\G$ through $(0,-k)=(0,0)$, namely the left 
%derivative $G'_{-}(0)$ that equals $b_0=V(0)$. The only difference is
%that for $k=0$ the supremum of~\eqref{slope} is not a maximum, unless
%$\G$ is a straight line close to the origin.

%facts established in the proof of Theorem~\ref{AJ} deserve
%emphasis: For $k\in(0,\kmax)$, the slope of the 
%supporting line through $(0,-k)$ to the curve $\G$ is equal to $V(k)$, and 
%$\theta_2<0$ maximises the fraction \eqref{slope} if and only if 
%$(\theta_2, G(\theta_2))$ is on that supporting line.
\end{remark}

\subsection{Almost worst case distributions}\label{sec-Almost worst case distributions}

For the problem \eqref{V}, call a density $p$ 
%Call a distribution $\Pd$ 
an {\em $(\epsilon$-$\gamma)$-Almost-Worst-Case-Density} (AWCD), where 
$\epsilon\ge 0,\;\gamma\ge 0$, if
%($\epsilon$-$\gamma$-AWCD) if 
\begin{equation}\label{AWC}
H( p )\leq k+\gamma\quad \mbox{and}\quad \int Xp d\mu \leq V(k) +\epsilon. 
\end{equation}
Thus, an $(\epsilon$-$\gamma)$-AWCD is a density which does not violate 
the constraint $H(p)\leq k$ by more than $\gamma$ and for which the expected
payoff does not exceed by more than $\epsilon$ the worst possible one subject
to the constraint. A {\em worst case density} (WCD) is a $(0$-$0)$-AWCD.

An $(\epsilon$-$\gamma)$ almost worst case distribution or a worst case
distribution is a distribution $\Pd$ 
whose density is an $(\epsilon$-$\gamma)$-AWCD or a WCD.

Theorem \ref{main} below establishes a clustering property of the 
$(\epsilon$-$\gamma)$-AWCDs, as well as a similar result for densities 
that almost attain the minimum in \eqref{W}.
 From a practical point of view, this may be
relevant for efficient hedging against the almost worst scenarios, but
this issue is not entered here.

Let us assign to each $\theta_2\in\Theta_2$ the unique
$\theta_1$ attaining the minimum in the 
definition~\eqref{L} of $G(\theta_2)$, determined in Lemma~\ref{lemma-L},  
% With this $\theta_1$, 
and denote
\begin{equation}\label{qte}
\qt(r):=p_\tete(r),\;\mbox{with $\theta_1$ attaining}\; 
K(\tete)-\theta_1=G(\theta_2).
\end{equation}
Given $k\in(0,\kmax)$, we will denote by $\qte$ the function $\qt$
with $\theta_2<0$ attaining the second maximum in Theorem~\ref{AJ}, i.e.,
%$\theta_1$ satisfying~\eqref{qte}. and $\theta_2<0$ attaining the second 
%This $\qte$ will be considered for $\theta_2<0$ attaining the 
%maximum in Theorem~\ref{AJ}. 
%By Remark~\ref{geom}, this maximising $\theta_2$ is
%the argument of that point of the graph of $G$ where 
%the supporting line through $(0,-k)$
%meets the graph (if there are several such points, either may be taken). 
%equivalently, the function %$\qte$ can be defined by
\begin{equation}\label{wcl}
\qte:=\qt=p_\tete \quad\mbox{with $(\tete)$ a maximiser in~\eqref{mm}}
\end{equation}
% i.e., with notation introduced before Lemma~\ref{lemma-T}, 
%for $\theta_2\in T(k)$. Then, 
%attains the double maximum in Theorem~\ref{AJ}. By Remark~\ref{geom},

\begin{theorem}\label{main}
(i) For $k\in(0,\kmax)$, %let $(\tete)$ attain the double maximum in 
%Theorem \ref{AJ}. Then
each $(\epsilon$-$\gamma)$-AWCD  $p$ belongs to the\\ Bregman neighborhood of
radius $(\gamma-\theta_2 \epsilon)$ of $\qte$ in~\eqref{wcl}, i.e., 
%in~\eqref{qte} with $\theta_2\in T(k)$, i.e., 
see \eqref{Bbe},     
\begin{equation}\label{Bbound}
B(p, \qte)\le \gamma-\theta_2 \epsilon\quad\mbox{if $p$ is an $(\epsilon$-$\gamma)$-AWCD}.
\end{equation}

(ii) For $\lambda >0$ with $-1/\lambda\in\Theta_2$, set 
%$\theta_1$ attain the mimimum in \eqref{L} for 
$\theta_2:=-1/\lambda$. % and $\theta_1$ satisfy~\eqref{qte}. 
Then for each density $p$
\begin{equation}\label{Wbound}
\int Xp d\mu +\lambda H(p)\ge W(\lambda)+\lambda B(p, \qt). 
\end{equation}
\end{theorem}
\begin {corollary}
\label{C.1}
Let $\{p_n\}$ be a sequence of $(\epsilon_n$-$\gamma_n)$-AWCDs with
\mbox{$\epsilon_n\to 0,\; \gamma_n\to 0$} in case (i), or a sequence of
densities with $\int Xp_n d\mu+\lambda H(p_n)\to W(\lambda)$ in case (ii).
Then $p_n$ converges to $\qte$ respectively to $\qt$ locally in 
measure.\footnote{This 
means that $\mu(\{r\in C: \vert p_n(r) - \qt (r) \vert > \epsilon
  \})\rightarrow 0$ for each $C\subset\Omega$ with $\mu( C )$ finite,
  and any $\epsilon > 0$. If $\mu$ is a finite measure, this is equivalent to 
standard (global) convergence in measure.} In particular, the function
$\qte$ %(or $p_\tete$ in case (ii)) 
is unique.
\end{corollary}
\begin{proof} (i) %The assumption implies $G(\theta_2)=K(\tete)-\theta_1$,
%see~\eqref{L}. Hence, using 
By the generalised Pythagorean identity 
Lemma~\ref{Pyti}, applied to $\tete$ in \eqref{wcl},
\begin{eqnarray}\label{Pbound}
H(p)\ge \theta_1+\theta_2\int Xpd\mu -K(\tete)+B(p,p_\tete) \nonumber    \\
   = \theta_2\int Xpd\mu-G(\theta_2)+B(p,\qt),
\end{eqnarray}
for each density $p$. %By the choice of $(\tete)$, 
As $\theta_2$ attains the maximum in Theorem~\ref{AJ}, here $\qt=\qte$ and
%$[k+G(\theta_2)]/\theta_2=V(k)$ thus 
$G(\theta_2)=\theta_2 V(k)-k$.
Hence, \eqref{Pbound} implies \eqref{Bbound} for each density $p$ which is  
an $(\epsilon$-$\gamma)$-AWCD, thus satisfies \eqref{AWC}.

(ii) In this case, \eqref{Pbound} holds as before.
%(this time not replacing $p_\tete$ by $\qte$). 
Multiplying it by
$\lambda=-1/\theta_2$ and using that $-\lambda G(-1/\lambda)=W(\lambda)$, see 
\eqref{W=}, we obtain \eqref{Wbound}.

The Corollary follows since $B(p_n, \qte)\to 0$ implies convergence of
$p_n$ to $\qt$ locally in measure \cite[Corollary 2.14]{CsiszarMatus2012art}. 
\end{proof}

\begin{remark}\label{Re2}
Corollary \ref{C.1} extends the known result that $\qte$ is a {\em generalized
solution} of problem \eqref{F} in the sense \cite{CsiszarMatus2012art}
that densities $p_n$ with
$\int Xp_n d\mu=b=V(k),\; H(p_n)\to F( b)=k$ converge to $\qte$ locally in 
measure,
and also establishes its (new) counterpart for problem \eqref{V}. 
\end{remark}

The function $\qte$ in Theorem~\ref{main}(i) will be called 
{\em worst case localiser}, for the almost worst case densities are clustering
in its (Bregman) neighborhood. This nice intuitive interpretation of
the function $\qte=p_\tete$ is complemented by the additional intuitive fact
that, by \eqref{Bbound}, its parameter $\theta_2$ controls the radius 
of that neighborhood. Most appealing is
the special case $\gamma=0$ of~\eqref{Bbound} that all densities that 
satisfy $H(p)\le k$ and yield expected payoff not excceding
the worst case by more than $\epsilon$, are contained in
a Bregman neighborhood of $\qte$ of radius
%\footnote{A smaller radius does not suffice, at least in those cases
%when $\theta_1+\theta_2 X(r)\ge \beta'(r,0)$ (thus the last term in
%\eqref{Py} vanishes), and some density satisfies 
%$H(p)= k,\:\int Xpd\mu=V(k)+\epsilon$.}
proportional to $\epsilon$, with proportionality factor $-\theta_2$.
The essence of the Corollary is that the Bregman distance of $p_n$
from $\qte$ goes to $0$. For certain choices of the integrand $\beta$ 
this implies convergence even in a stronger sense than locally in measure,
see Example~\ref{Kullback}.   

Clearly, the worst case localiser $\qte$
coincides with the WCD whenever the latter exists (apply 
\eqref{Bbound} to $\epsilon=\gamma=0$). A necessary and sufficient condition 
for the existence of a WCD is given in Lemma~\ref{l-opcond}(ii). There, we
have skipped the proof that 
a WCD has to be of form $p_\tete$ with $\theta_2<0$, which is obvious now
as the WCD is a worst case localiser. Lemma~\ref{l-noc} below will also be
useful in identifying situations when the worst case localiser is actually a
WCD. Its Corollary addresses the simplest such situation. 
% which, in effect, has been covered already in~\cite[Theorem 1]{BreuerCsiszar2013}.
%$\theta_2<0$. %appears in \cite[Theorem
%1]{BreuerCsiszar2013-MaFi}; this issue will be revisited in 
%see Subsection~\ref{largek}.
When in \eqref{V} the minimum is not attained,
the worst case localiser may or may not be a density, see the examples
below, though 
it always satisfies $\int \qte d\mu \le 1$, see Theorem~\ref{AJ}. 
Note that the computation of the worst case localiser is not harder than
the computation of $V(k)$ along the lines of Subsection~\ref{sec-Ahmadi}, 
for that calculation does provide the parameters $\tete$ of
$\qte=p_\tete$ that attain the double maximum in Theorem~\ref{AJ}.

\begin{lemma}\label{l-noc}
A function $p_\tete$ in~\eqref{EF} with %$(\tete)\in \Theta,\;
$\theta_2<0$ is the worst case
localiser $\qte$ for $k\in (0,\kmax)$ if and only if the vector
$(1-\int p_\tete d\mu,\, V(k)-\int Xp_\tete d\mu)\in\R^2$ belongs to the
normal cone of $\dom\, K$ at $(\tete)$, i.e., for each $(\tebteb)\in\dom\, K$
\begin{equation}\label{noc}
(\bar{\theta}_1-\theta_1)\left(1-\int p_\tete d\mu\right)+(\bar{\theta}_2-\theta_2)
\left(V(k)-\int Xp_\tete d\mu\right) \le 0.
\end{equation}
\end{lemma}
\begin{corollary} If the worst case localiser $\qte=p_\tete$ has parameters 
$(\tete)$  in the interior of $\dom\,K$, then it is a WCD.
%$(\tete)\in\intdom K$ then it is a WCD.
%A sufficient condition for the worst case localiser
%$p_\tete$ to be a WCG is $(\tete)\in\intdom K$.
\end{corollary}
\begin{proof} By Theorem~\ref{AJ}, the condition in the definition~\eqref{wcl}
%The worst case localiser is defined by the condition that
%$(\tete)$ attains the maximum in~\eqref{mm} or
is equivalent to the condition that $(\tete)$ attains 
the maximum of $f(\tete):=\theta_1+\theta_2 b-K(\tete)$ where
$b=V(k)$. The latter is satisfied if and only if for each 
$(\tebteb)\in\dom\,K$, the concave function
%$$f(t):=\theta_1+t(\bar{\theta}_1-\theta_1)+[\theta_2+t(\bar{\theta}_2-
%\theta_2)]b-K(\theta_1+t(\bar{\theta}_1-\theta_1),\theta_2+t(\bar{\theta}_2-
$$ f(t):=f(\theta_1+t(\bar{\theta}_1-\theta_1),\;\theta_2+t(\bar{\theta}_2-
 \theta_2)),\quad 0\le t\le 1 $$
is maximised by $t=0$, i.e., its (right) derivative at $t=0$
is nonpositive. On account of~\eqref{dide}, that condition is equivalent to 
\eqref{noc}. The Corollary follows since the normal cone of $\dom\, K$
at an interior point consists of $(0,0)$ alone. Thus Lemma~\ref{l-noc}
gives the conditions $\int p_\tete d\mu=1, \int Xp_\tete d\mu=V(k)$, which
mean by Lemma~\ref{l-opcond} that $p_\tete$ is a WCD.
\end{proof}

\begin{example}\label{Kullback} Let $H(p)$ be the $I$-divergence,
%$D(\Pd||\Pd_0)$, 
formally (see Example 1) let $\beta$ be the autonomous
integrand given by $f(s)=s\log s$ and let $\Pd_0=\mu$. Then $f^*(\tau)=e^{\tau-1}$,
$K(\tete)=\int e^{\theta_1+\theta_2 X-1}d\mu$, and
$$G(\theta_2)=\min_{\theta_1}[K(\tete)-\theta_1]=\log \int e^{\theta_2X} d\mu
:=\Lambda(\theta_2),$$ 
for $\theta_2\in\Theta_2=\dom\:\Lambda$. The minimum in the definition of
$G(\theta_2)$ is attained for $\theta_1=1-\Lambda(\theta_2)$. If 
$m>-\infty$ and 
 $X(r)=m$ on a set of $\mu$-measure $\mu_0>0$ then\footnote{Note that
Theorems~\ref{AJ},~\ref{main} do not apply to $k=\kmax$. Here, in that case  
the WCD equals $1/\mu_0$ on the set $\{r:X(r)=m\}$ 
and $0$ elsewhere. It does not belong to the family~\eqref{EF}, and the almost 
worst case densities do not cluster in its Bregman neighborhood.}
$\kmax=-\log \mu_0$, 
otherwise $\kmax=+\infty$ (assuming \eqref{nont}). %i.e., that $\dom\,\Lambda$
%contains some $\theta_2<0$). 

For $k\in (0,\kmax)$,
Theorem~\ref{AJ} gives
$V(k)=\max_{\theta_2<0}[k+\Lambda(\theta_2)]/\theta_2$ 
and Theorem~\ref{main} gives the worst case
localiser $\qte=\exp (-\Lambda(\theta_2)+\theta_2X)$, with 
$\theta_2$ attaining the
above maximum. %(the parameter $\theta_1=1-\Lambda(\theta_2)$ is suppressed). 
This worst case localiser is always a density. It also
satisfies $H(\qte)=k$ and hence is actually a WCD, 
except for the case when $\dom\,\Lambda$ contains its left endpoint 
$\thmin$, $\Lambda'(\thmin)$ is finite, and
$k>H(q_{\thmin})$. %As pointed out before 
%\cite{BreuerCsiszar2013-MaFi}, in that exceptional case the minimum 
%in~\eqref{V} is not attained. 
In that case the maximiser in
Theorem~\ref{AJ}, equal to the parameter of the worst case 
localiser, is $\theta_2=\thmin$. Note that for this example,  
the formula for $V(k)$ appears in
\cite{Ahmadi-Javid2011} and \cite{BreuerCsiszar2013-MaFi} show that in the above exceptional 
case the minimum in \eqref{V} is not attained. The result of 
Theorem~\ref{main} appears new even for this special case.

In this example, the Bregman distance 
\eqref{Bbe} of densities coincides with $I$-divergence, hence 
Theorem~\ref{main} gives $D(p||\qte)\le\gamma-\theta_2 \epsilon$
for each $(\epsilon$-$\gamma)$-AWCD $p$.
%satisfies $D(p||p_{\theta_2})\le\gamma-\theta_2 \epsilon$. 
In the Corollary of Theorem~\ref{main}, now the  
almost worst case densities converge to the worst case localiser
in a much stronger sense than in measure. Indeed, the result that 
their $I$-divergence from the worst case localiser approaches $0$ 
is stronger than $L_1(\mu)$ convergence to the worst case localiser.
\end{example}

\begin{example}
\label{Burg} Again in the setting of Example 1, take now
$f(s)=-\log s$. Then $H(p)$ equals reverse $I$-divergence, i.e., the 
$I$-divergence of the default distribution $\Pd_0=\mu$ from the 
distribution $\Pd$ with density $p$. As $f$ is not cofinite, the
standing assumption $\kmax>0$ holds if and only if $m>-\infty$.

Take specifically $\Omega=(0,1),\;X(r)=r$, and take for $\mu=\Pd_0$ the
distribution with Lebesgue density $2r$. As
$f^*(\tau)=-1-\log(-\tau)\;(\tau<0)$, then $K(\tete)=\int_0^1
[-1-\log(-\theta_1-\theta_2 r)]\,2rdr$ and $p_\tete(r)=1/(-\theta_1-\theta_2 r)$ 
for $(\tete)\in\Theta=\dom\,K=\{(\tete):\theta_1\le
0,\,\theta_1+\theta_2<0\}$. Simple calculus shows that for 
$-2\le\theta_2\le 0$ the minimum in the 
definition of $G(\theta_2)$ is attained for $\theta_1$ such that
$p_\tete=\qt$ is a density, but the functions $G$ and $\qt$
can not be given explicitly.
If $\theta_2<-2$
then this minimum is attained for $\theta_1=0$, and 
$G(\theta_2)=-\frac{1}{2}-\log (-\theta_2)$.  
%Assigning to $\theta_2$ the minimising $\theta_1$, 
One sees that
$H(\qt)$ ranges from $0$ to $\log 2-1/2$ as $\theta_2$ ranges from $0$ to
$-2$. Hence in case $k\le \log 2-1/2$ the WCD exists, it equals 
that $p_\tete$ which is a density and satisfies $H(p_\tete)=k$.
In case $k\ge\log2-1/2$ the worst case localiser is %the function
$\qte(r)=-1/\theta_2 r$  with $\theta_2\le -2$ %determined from the 
%condition in Theorem 1 that it 
attaining $V(k)=\max_{\theta_2<0}[k+G(\theta_2)]/\theta_2$. By
simple calculus, this maximiser is $\theta_2=-e^{k+1/2}$, the maximum is
$V(k)=e^{-(k+1/2)}$, and the worst case localiser is %the function
$\qte(r)=\frac{1}{r}e^{-(k+1/2)}$, which is not a 
density unless $k= \log2-1/2$.  

In this case the Bregman distance~\eqref{Bbe} is
$$ B(p,q)=\int \left[\log \frac{q}{p}+\frac{p}{q}-1\right] d\mu.$$
\noindent The Corollary of Theorem~\ref{main} now
does not admit a substantial strengthening, for the result that
this Bregman distance approaches $0$ does not imply 
convergence in a familiar sense stronger than in measure.
\end{example}

\begin{example}\label{never} Let $\Omega,\; X,\; \mu$ be as in 
Example~\ref{Burg}, but this time let the default distribution $\Pd_0$ 
be the uniform distribution whose $\mu$-density is $p_0(r)=\frac{1}{2r}$.
Take for $H(p)$ the Bregman distance $B(p,p_0)$ in Example~\ref{Burg},
i.e., the integral functional~\eqref{H1}
with $\beta(r,s)=\Delta_f(s,p_0(r))=-\log s-\log(2r)+2r(s-\frac{1}{2r})$.
Then $\beta^*(r,\tau)=\log{2r}-\log(-\tau+2r)$,
$(\beta^*)'(r,s)=1/(-\tau+2r)$, $\tau<2r$.
 The set $\Theta$
(equal to $\dom\, K$) of this example consists of those $(\tete)$ for which
$(\theta_1, \theta_2-\theta_1)$ belongs to the set $\Theta=\dom\,K$ of 
Example~\ref{Burg}. Moreover, for such $(\tete)$ the function 
$p_\tete(r)= 1/[-\theta_1-(\theta_2-2)r]$ coincides with the function
$p_{\theta_1,\theta_2-2}$ of Example~\ref{Burg}, which
 can not be a density if $\theta_2<0$. This proves that
in the present Example no WCD exists for any $k>0$.
\end{example}

\subsection{Effect of the threshold $k$ on the existence of a WCD}
\label{sec-largek}

%Given the integrand $\beta$ and the payoff function $X$, it
%may depend on the choice of the threshold $k$  whether a worst case
%density attaining the minimum in \eqref{V} exists (in which case 
%worst case localiser equals the WCD). Below, $k\in (0,\kmax)$ is assumed,
%even if not stated explicitly.
\rm
This subsection addresses the effect of the choice of the threshold $k$ on the 
worst case localiser, in particular on whether that localiser is also a WCD.
%Below, $k\in (0,\kmax)$ is assumed, even if not stated explicitly.

Examples \ref{Kullback}, \ref{Burg}, and \ref{never} demonstrate that a WCD may exist for all or for no $k$, or
there may exist a critical value $\kcr$ such that a WCD exists if 
$k<\kcr$ but does not exist if $k>\kcr$. It appears a plausible conjecture
that these three alternatives are exhaustive, i.e., that if
a WCD exists for some $k$, it also exists for each $k'<k$. While this
conjecture remains open in general, it will be proved under 
conditions that cover many typical cases. 

%Below, $k\in (0,\kmax)$ is assumed, even if not stated explicitly. 
Recall that $\thmin$ with $-\infty\le \thmin<0$ denotes the left endpoint of
the interval $\Theta_2$, the projection of $\dom\, K$ to the $\theta_2$ axis.
The condition $m>-\infty$ is necessary for $\kmax<+\infty$ and 
sufficient for $\thmin=-\infty$.

%Part (ii) of the next Theorem gives a general sufficient condition for the 
%worst case localiser to behave as in Example~\ref{Kullback}.

\begin{theorem}\label{WCLWCD} (i) If for some $k\in(0,\kmax)$ the worst case
localiser $\qte$ is a density, it is a WCD for $k$  
unless\footnote{Here $G'(\thmin)$ means the right derivative.} %$G'_{+}(\thmin)$.}
\begin{equation}\label{derfin}
\thmin\in\Theta_2,\quad G'(\thmin)>-\infty
\end{equation}
and 
\begin{equation}\label{kcr}
k>\kcr:= -G(\thmin)+\thmin G'(\thmin).
\end{equation}
If~\eqref{derfin} and~\eqref{kcr} hold then $\qte=\hat{q}_{\kcr}=q_{\thmin}$
and no WCD exists for $k$.

(ii) If $\dom\,K$ contains the $\theta_1$-axis, i.e., 
$\int \beta^*(r,\theta_1)\mu(dr)$ is finite for each $\theta_1\in\R$, then
the worst case localiser $\qte$ is a density and hence it is a WCD unless \eqref{derfin} and \eqref{kcr} hold.
\end{theorem}
\begin{proof} (i) Suppose $\qte=p_\tete$ in~\eqref{wcl} is a density. By
Lemma~\ref{l-opcond} it is a WCD for $k$ if and only if
\begin{equation}\label{XV}
\int Xp_\tete d\mu =V(k).
\end{equation}
Since $ p_\tete$ is a worst case localiser and $\int p_\tete=1$, Lemma~\ref{l-noc}
gives that
\begin{equation}\label{rnoc}
(\bar{\theta}_2-\theta_2) \left(V(k)-\int X p_\tete d\mu\right)\le
0\quad\mbox{for}\quad\bar{\theta}_2\in \Theta_2.
\end{equation}
This immediately implies ~\eqref{XV} if $\theta_2\neq \thmin$, or equivalently
(see Remark~\ref{geom}) if the supporting line through $(0,-k)$ to the curve
$\G$ does not contain $(\thmin, G(\thmin))$. This is always the case if 
\eqref{derfin} does not hold, and also when~\eqref{derfin} holds but
$G'(\thmin)$, the largest slope of supporting lines to $\G$ at
$(\thmin, G(\thmin))$, is less than $[k+G(\thmin)]/\thmin$. As the last
condition is equivalent to $k<\kcr$, only the case $k=\kcr$ remains to
cover to complete the proof
that $\qte$ is a WCD unless~\eqref{derfin} and~\eqref{kcr} hold.

In that remaining case, $\qte=p_\tete$ with $\theta_2=\thmin$, and
instead of~\eqref{XV} only the inequality $\int X\ptt d\mu \ge V(\kcr)$ 
follows from~\eqref{rnoc}. Suppose indirectly that it
is strict, then Lemma~\ref{Lemma1} implies that 
$\int X\ptt d\mu=V(\bar{k})$ for 
some $\bar{k}\in(0,\kcr)$ (as the integral is less than $b_0$ by 
Lemma~\ref{l-opcond}). This means by Lemma~\ref{l-opcond} that
 $p_\tete$ (with $\theta_2=\thmin$) is a WCD for 
$\bar{k}$. Hence, by Remark~\ref{geom}, the supporting line through
$(0,-\bar{k})$ to the curve $\G$ contains $(\thmin, G(\thmin))$, 
contradicting the fact that among the supporting lines to
$\G$ at $(\thmin, G(\thmin))$ the one through $(0,-\kcr)$ has the largest
slope. This contradiction proves that~\eqref{XV} holds and hence $\qte$
is a WCD also when $k=\kcr$.

The last assertions of part (i) are obvious. Indeed, if~\eqref{derfin}
holds and $k>\kcr$ then the supporting line through $(0,-k)$ to $\G$
meets the curve at $(\thmin, G(\thmin))$, just as the supporting line
through $(0,-\kcr)$ does, hence $\qte=q_{\thmin}=q_{\kcr}$. As it is a WCD for
$\kcr$, it can not be a WCD for $k>\kcr$ with $V(k)<V(\kcr)$.

(ii) The hypothesis implies that a vector in $\R^2$ can belong to the
normal cone of $\dom\, K$ at some $(\tete)$ only if the first component of
this vector is $0$. On account of Lemma~\ref{l-noc}, this proves that 
the worst case localiser $\qte=p_\tete$, for any $k\in(0,\kmax)$, has to
satisfy $\int p_\tete d\mu =1$.
\end{proof}
\begin{corollary}\label{diff} The function $G$ is differentiable at each
  $\theta_2\in(\thmin, 0)$ for which $\qt$ in~\eqref{qte} is a density.
\end{corollary}
\begin{proof} Suppose indirectly that the curve $\G$ has several supporting
lines at $(\theta_2, G(\theta_2))$, say one containing $(0,-k_1)$ and
another $(0,-k_2)$, where $k_1\neq k_2$. Then
$\hat{q}_{k_1}=\hat{q}_{k_2}=\qt$, see Remark~\ref{geom}, hence $\qt$ is the
WCD both for $k_1$ and $k_2$, by Theorem~\ref{WCLWCD}. This means that
$V(k_1)=\int X\qt d\mu =V(k_2)$, contradicting  $k_1\neq k_2$.
\end{proof}

Finally, we discuss for $f$-divergence balls, see Example 1, the 
dependence on the threshold $k$ (the ``radius'' of the ball)
of the worst case localiser and whether it is a WCD. Formally, let
$\beta(r,s)=f(s)$ be an autonomous integrand, $f$ strictly
convex and differentiable on $(0,+\infty)$, $f(0)=\lim_{s\downarrow 0} f(s)$,
$f(1)=0$, let $\mu$ be a probability measure, and $\Pd_0=\mu$.

The case of cofinite $f$ is covered by Theorem~\ref{WCLWCD}(ii), the
integral in its hypothesis being equal to $f^*(\theta_1)$, finite for each 
$\theta_1\in\R$. Therefore we focus on the non-cofinite case, supposing
\begin{equation}\label{c}
\lim_{s\uparrow +\infty}\frac{f(s)}{s}=c,\quad\quad c\quad\mbox{finite}.
\end{equation}
Then the standing assumption $\kmax>0$ (equivalent to $\thmin<0$)
holds if and only if $m>-\infty$. With no loss of
generality, assume that $m=0$ (clearly, the minimisation problem~\eqref{V}
is not affected by adding a constant to $X$).

Under the above assumptions, $K(\tete)=\int f^*(\theta_1+\theta_2 X) d\mu$
with $\theta_2<0$ is finite if $\theta_1<c$ and infinite if $\theta_1>c$,
because~\eqref{c} implies that $f^*(\tau)$ is finite for $\tau<c$ but not for
$\tau>c$. %($(c,\theta_2)$ may or may not belong to $\dom\, K$).
It follows for any $\theta_2<0$ that the associated $\tilde{\theta}_1$ in Lemma
\ref{lemma-L} is equal to $c$, hence Lemma~\ref{lemma-L} gives that 
the function $\qt=p_\tete$ in~\eqref{qte} is a density if and only if
\begin{equation}\label{integ}
g(\theta_2):=\int (f^*)'(c+\theta_2 X) d\mu \ge 1.
\end{equation}

Moreover, if $g(\theta_2)\le 1$ then 
\begin{equation}
\label{eq-qt2}
\qt =(f^*)'(c+\theta_2 X)=p_{c,\theta_2}.
\end{equation}
%Define\footnote{If $q_{\tilde{\theta}_{\min}}$ is a density then 
%$G'(\tilde{\theta}_{\min})$ in~\eqref{crtilde} is well defined
%by Corollary~\ref{diff}; otherwise it is interpreted as right 
%derivative.} 

%Denote by $\tilde{\theta}_{\min}$ the infimum of those $\theta_2<0$ for
%which~\eqref{integ} holds, and if $\tilde{\theta}_{\min}>-\infty$, set
%\footnote{If $G$ is not differentible at $\tilde{\theta}_{\min}$, 
% take the right derivative in \eqref{crtilde}.} 
%$G'(\tilde{\theta}_{\min})$ in~\eqref{crtilde} is well defined
%by Corollary~\ref{diff}; otherwise it is interpreted as right 
%derivative
%\begin{equation}\label{mintilde}
%\tilde{\theta}_{\min}:=\inf\{\theta_2:\int(f^*)'(c+\theta_2 X) d\mu \ge 1\} 
%\end{equation}
%\begin{equation}\label{crtilde}
%\tilde{k}_{\mathrm{cr}}:=\tilde{\theta}_{\min} G'_{+}(\tilde{\theta}_{\min})-
%G(\tilde{\theta}_{\min}).
%\end{equation}
%This has to be distinguished from $\kcr$ in
%Theorem~\ref{WCLWCD} (now undefined) but will play a similar role.
\begin{theorem}\label{fballs}
Under the assumptions~\eqref{c} and $m=0$, if %for each $\theta_2<0$
%\begin{equation}\label{infini}
%\int (f^*)'(c+\theta_2 X) d\mu=+\infty
%\end{equation}
%the integral in~\eqref{integ} equals 
$g(\theta_2)=+\infty$ for each $\theta_2<0$ 
then the WCD exists for all $k\in (0,\kmax)$. Otherwise, denote
\begin{equation}\label{tetil}
\tetil:=\inf\{\theta_2 : g(\theta_2)\ge 1\},
\end{equation}
\begin{equation}\label{crtilde}
\tilde{k}_{\mathrm{cr}}:=\tilde{\theta}_{\min} G'_{+}(\tilde{\theta}_{\min})-
G(\tilde{\theta}_{\min}).
\end{equation}
Then $\tetil\in (-\infty,0)$, $\tilde{k}_{\mathrm{cr}} \in (0,\kmax)$, and
for $k<\tilde{k}_{\mathrm{cr}}$ the WCD exists. For 
$k\ge\tilde{k}_{\mathrm{cr}}$ the worst case
localiser $\qte$ is of form \eqref{eq-qt2}; it is not a density if
$k>\tilde{k}_{\mathrm{cr}}$, while for $k=\tilde{k}_{\mathrm{cr}}$ it is a
density (and hence a WCD) unless\footnote{It is
left open whether that exceptional case is possible.} $g(\theta_2)<1$   
for each $\theta_2<0$ with $g(\theta_2)<+\infty$.
\end{theorem}
\begin{proof} Since in the current case $\thmin=-\infty$, 
Theorem~\ref{WCLWCD} implies that
the worst case localiser $\qte$ is a WCD if and only if it is a density.
By the passage preceding the Theorem, the latter holds if 
and only if $\qte=\qt$ with $\theta_2$
satisfying~\eqref{integ}. This immediately proves the first assertion.

Suppose next that $g(\theta_2)$ is finite for some $\theta_2<0$, and
denote the %integral in~\eqref{integ} by $g(\theta_2)$ and the
supremum of such parameters $\theta_2$  by $\sigma$.
One verifies via monotone convergence and dominated convergence that
$g(\theta_2)$  is a continuous, strictly increasing function of $\theta_2\in
(-\infty, \sigma)$ that approaches $0$ or $g(\sigma)$
as $\theta_2$ goes to $-\infty$ or $\sigma$. Hence if 
$g(\sigma)\ge 1$ then $\tilde{\theta}_{\min}$ is equal to the 
unique $\theta_2\le\sigma$ with $g(\theta_2)=1$, whereas
if $g(\sigma)<1$ then $\tilde{\theta}_{\min}=\sigma$. 
In both cases $-\infty<\tilde{\theta}_{\min}<0$,
using that $g(0)=+\infty$. 

By the definition~\eqref{crtilde} of $\ktil$, the supporting line
to $\G$ at 
$(\tilde{\theta}_{\min}, G(\tilde{\theta}_{\min}))$ 
of slope $G'_{+}(\tilde{\theta}_{\min})$ intersects the vertical
axis at $(0,-\tilde{k}_{\mathrm{cr}})$. Hence
$\tilde{k}_{\mathrm{cr}}>0$, unless the function $G$ is linear in the interval
$[\tilde{\theta}_{\min} ,0]$; the latter possibilty will be ruled out in 
the Appendix. It follows, too, that the supporting line to $\G$
through $(0,-k)$ with $k<\tilde{k}_{\mathrm{cr}}$ or 
$k>\tilde{k}_{\mathrm{cr}}$
meets the curve $\G$ at a point (or points) with argument 
$\theta_2>\tilde{\theta}_{\min}$ respectively 
$\theta_2\le\tilde{\theta}_{\min}$. Moreover, the latter inequality is
strict if $g(\tilde{\theta}_{\min})=1$ (equivalent to $g(\sigma)\ge 1$), for 
in that case 
$G$ is differentiable at $\tilde{\theta}_{\min}$ due to Corollary~\ref{diff}.

Referring to Remark~\ref{geom}, the above considerations prove that 
the parameter $\theta_2$ in the 
representation $\qte=\qt$ in~\eqref{wcl} satisfies or does not satisfy the
condition~\eqref{integ} if $k<\tilde{k}_{\mathrm{cr}}$ respectively
$k>\tilde{k}_{\mathrm{cr}}$, no matter whether
$g(\sigma)\ge 1$ or not. These facts, and that $\qte=q_{\tilde{\theta}_{\min}}$
if $k=\tilde{k}_{\mathrm{cr}}$, imply all remaining assertions of the 
Theorem, see the first passage of the proof.
\end{proof}

% \newpage
\section*{Appendix}

\subsection*{Proof of Lemma~\ref{lemma-L}} 
\begin{proof} 
Fix $\theta_2\in\Theta_2$, define
$\tilde{\theta}_1$ as in the lemma. Then
$(\tete)\in\Theta$ for all
$\theta_1<\tilde{\theta}_1$, and the function $f(\theta_1):=K(\tete)$ is
convex, closed, and differentiable in its effective domain 
$(-\infty,\tilde{\theta}_1 )$, with
\begin{equation}\label{par3}
f'(\theta_1)=\int p_\tete d\mu,\quad \theta_1<\tilde{\theta}_1,
\end{equation}
see~\eqref{par1}. If $(\tilde{\theta}_1,\theta_2)\in\Theta$ then \eqref{par3}
holds also for the left derivative at $\theta_1=\tilde{\theta}_1$. Hence the
last assertion of the Lemma immediately follows.

To prove that one of the alternatives (i) and (ii) indeed takes place, note
that the properties of $\beta^*$ stated in the passage after \eqref{bstar}
imply, by monotone convergence, that $f'(\theta_1)$ in \eqref{par3} goes to
$0$ if $\theta_1\downarrow -\infty$ and to $+\infty$ if
$\tilde{\theta}_1=+\infty$ and $\theta_1\uparrow +\infty$. Hence, due to
continuity of $f'(\theta_1)$, alternative (i) fails only if
\begin{equation}\label{le1}
\int p_\tete d\mu<1\;\mbox{for all}\;\theta_1\;\mbox{with}\;(\tete)\in\Theta,
\end{equation}
and \eqref{le1} can hold only if $\tilde{\theta}_1<+\infty$. Further,
\eqref{le1} implies that $(\tette)\in \dom\;K$, for in the opposite case 
$f(\tilde{\theta}_1)=+$ the derivative \eqref{par3} of the closed convex
function $f(\theta_1)$ would go to $+\infty$ as
$\theta_1\uparrow\tilde{\theta}_1$.

The proof will be complete if we show that \eqref{le1} implies 
$(\tette)\in\Theta_2$. It has already been shown to imply 
$(\tette)\in \dom\;K$, in particular,
that $\tilde{\theta}_1+\theta_2 X(r)\le\beta'(r,+\infty)\;\mu$-a.e., 
thus it remains to verify, see \eqref{Th}, that the set
$\{r:\tilde{\theta}_1+\theta_2 X(r)=\beta'(r,+\infty)\} $
has $\mu$-measure $0$. On that set,
$p_\tete(r)=(\beta^*)'(r, \theta_1+\theta_2 X(r))$ grows to $+\infty$ as
$\theta_1\uparrow \tilde{\theta}_1$. Hence, were it not a $0$-measure set,
$\int p_\tete d\mu$ would grow to $+\infty$, contradicting \eqref{le1}.
\end{proof}

\subsection*{Proof of Lemma~\ref{Fstar}} 
\begin{proof} 
Fix $\theta_2 \in\R$ and consider the 
(not necessarily proper) convex function 
$$ L(a) := \inf_{b\in\R} (J(a,b)-\theta_2 b),\quad a\in\R.$$
Then
\begin{eqnarray*} 
F^*(\theta_2)=\sup_b (\theta_2 b-F(b))=-\inf_b(F(b)-\theta_2b)=-L(1)\\
 =-L^{**}(1)=-\sup_{\theta_1}(\theta_1-L^*(\theta_1))=
    \inf_{\theta_1}(L^*(\theta_1)-\theta_1),
\end{eqnarray*}
where the third equality holds since $F(b)=J(1,b)$, and the 
fourth one holds since $a=1$ is in the interior of $\dom\;L$.
Here
\begin{eqnarray*}
L^*(\theta_1)=\sup_a (\theta_1 a-L(a))=\sup_a[\theta_1 a+\sup_b
(-J(a,b)+\theta_2 b)]\\=J^*(\tete)=K(\tete).
\end{eqnarray*}
Recalling the definition \eqref{L} of $G$, this completes the proof.
\end{proof}

\subsubsection*{Completion of the proof of Theorem~\ref{fballs}}
It remains to rule out the
possibility that the function $G$ is linear in the interval $[\tetil,0]$.
Suppose indirectly that for some $\tilde{b}\in\R$
\begin{equation}\label{indi}
G(\theta_2)=\theta_2 \tilde{b}\quad\mbox{if}\quad \theta_2\in [\tetil,0].
\end{equation}
Here necessarily $\tilde{b}\le b_0$, by~\eqref{Gbound}. As~\eqref{indi} 
implies $G^*(\tilde{b})=0$, which means by~\eqref{FL} that $F(\tilde{b})=0$,
it follows by Remark 2 that actually $\tilde{b}=b_0$.

As $g(\tetil)\le 1$ by the proof of Theorem~\ref{fballs}, the value
$\theta_1$ in~\eqref{qte} attaining $K(\tete)-\theta_1=G(\theta_2)$ for 
$\theta_2=\tetil$ is equal to $c$. Thus
Lemma~\ref{Pyti} applied to $p=p_0$ and $(\tete)=(c,\tetil)$ gives
$$ 0=H(p_0)\ge \tetil\int Xp_0 d\mu-G(\tetil)+B(p_0, q_{\tetil}).    $$
Here the integral equals $b_0$ by definition, and $\tilde{b}$ in~\eqref{indi}
has been shown to equal $b_0$. Hence it follows that $B(p_0, q_{\tetil})=0$, 
which means that $q_{\tilde{\theta}_{\min}}$ equals $p_0=1$ ($\mu$-a.e.).
By Remark \ref{unique}, this contradicts $\tilde{\theta}_{\min}\neq 0$, proving that the indirect assumption~\eqref{indi} is false.

\bibliographystyle{plain}
\bibliography{GenEnt}

\end{document}